\documentclass[11pt]{article}
\usepackage{mathrsfs}
\usepackage{amsfonts}

\usepackage{amsfonts, amsmath, amssymb}
\usepackage{amssymb,amsfonts,amsmath,
latexsym, epsfig,cite, psfrag,eepic,colordvi}
\usepackage{amscd,graphics}
\usepackage{lineno}

\newcommand{\qed}{\hfill $\Box $}
\newcommand{\pf}{\noindent {\bf Proof.} }

\newtheorem{theorem}{Theorem}[section]
\newtheorem{lemma}[theorem]{Lemma}
\newtheorem{coro}[theorem]{Corollary}

\newtheorem{conjecture}[theorem]{Conjecture}

\newtheorem{problem}{Problem}

\begin{document}

\title{Vertex-Coloring Edge-Weighting of Bipartite Graphs  with Two
Edge Weights
\thanks{This work was supported in part by  the National Natural
Science Foundation of China No. 11101329.}}
\author{Hongliang Lu\footnote{E-mail address:
luhongliang215@sina.com}
\\ {\small School of Mathematics and Statistics}
\\ {\small Xi'an Jiaotong University, Xi'an, 710049, PR China}
}

\date{}

\maketitle

\begin{abstract}
Let $G$ be a graph and  $\mathcal {S}$ be a subset of $Z$. A
vertex-coloring $\mathcal {S}$-edge-weighting of $G$ is an
assignment of weight $s$ by the elements of $\mathcal {S}$ to each
edge of $G$ so that adjacent vertices have different sums of
incident edges weights.

It was proved that every 3-connected bipartite graph admits a
vertex-coloring $\{1,2\}$-edge-weighting (Lu, Yu and Zhang, (2011)
\cite{LYZ}). 
In this paper, we show that the following result: if a
3-edge-connected bipartite graph $G$ with minimum degree $\delta$
contains a vertex $u\in V(G)$ such that $d_G(u)=\delta$ and $G-u$ is
connected, then $G$ admits a vertex-coloring $\mathcal
{S}$-edge-weighting for $\mathcal {S}\in \{\{0,1\},\{1,2\}\}$. In
particular, we show that every 2-connected and 3-edge-connected
bipartite graph admits a vertex-coloring $\mathcal
{S}$-edge-weighting for $\mathcal {S}\in \{\{0,1\},\{1,2\}\}$. The
bound   is sharp, since there exists a family of infinite bipartite
graphs which are 2-connected and do not admit  vertex-coloring
$\{1,2\}$-edge-weightings or vertex-coloring
$\{0,1\}$-edge-weightings.
%
\\
[2mm] \textbf{Keywords}: edge-weighting; vertex-coloring; $2$-connected; bipartite graph.\\
[2mm] {\bf AMS subject classification (2000)}: 05C15.
\end{abstract}

\section{Introduction}

In this paper, we consider only finite, undirected and simple
connected graphs. For a vertex $v$ of graph $G=(V,E)$, $N_G(v)$
denotes the set of vertices which are adjacent to $v$ and
$d_G(v)=|N_G(v)|$ is called the \emph{degree} of vertex $v$. Let
$\delta(G)$ and $\triangle(G)$ denote the minimum degree and maximum
degree of graph $G$, respectively.
 For $v\in V(G)$ and $r\in Z^+$, let
$N_G^{r}(v)=\{u\in N(v)\ |\ d_G(u)=r\}$.  If $v\in V(G)$ and $e\in
E(G)$, we use $v\sim e$ to denote that $v$ is an end-vertex of $e$.
For two disjoint subsets $S,T$ of $V(G)$, let $E_G(S,T)$ denote the
subset of edges of $E(G)$ with one end in $S$ and other end in $T$
and let $e_G(S,T)=|E_G(S,T)|$.

Let   $\mathcal {S}$ be a subset of $Z$. A \emph{$\mathcal
{S}$-edge-weighting} of graph $G$ is an assignment of weight $s$ by
the elements of $\mathcal {S}$ to each edge of $G$.  A $\mathcal
{S}$-edge-weighting $w$ of   a graph $G$ induces a coloring of the
vertices of $G$, where the color of vertex $v$, denoted by $c(v)$,
is $\sum_{e\sim v}w(e)$. A $\mathcal {S}$-edge-weighting of a graph
$G$ is a \emph{vertex-coloring}   if for every edge $e = uv$,
$c(u)\neq c(v)$ and  we say that $G$ admits a \emph{vertex-coloring
$\mathcal {S}$-edge-weighting}.
  If $\mathcal {S}=\{1,2,\ldots,k\}$, then a
vertex-coloring $\mathcal {S}$-edge-weighting of graph $G$ is
usually called a \emph{vertex-coloring $k$-edge-weighting}.


For vertex-coloring edge-weighting, Karo\'nski, Luczak and Thomason
\cite{KLT} posed the following conjecture:
\begin{conjecture}\label{conj1}
Every graph with no isolated edge admits a vertex-coloring
$3$-edge-weighting.
\end{conjecture}
This conjecture is still wide open. Karo\'nski et al. \cite{KLT}
showed that the Conjecture \ref{conj1} is true for 3-colorable
graphs. Recently, Kalkowski et al. \cite{KKF} showed that every
graph with no isolated edge admits a vertex-coloring
5-edge-weighting. This result is an improvement on the previous
bounds on $k$ established by Addario-Berry et al. \cite{ADMRT},
Addario-Berry et al. \cite{ADR}, and Wang et al. \cite{WY}, who
obtained the bounds $k = 30, k = 16$, and $k = 13$, respectively.

Many graphs actually admit a vertex-coloring 2-edge-weighting (in
fact, experiments suggest (see \cite{ADR}) that almost all graphs
admit  a vertex-coloring 2-edge-weighting), however it is not known
which ones are not.  Chang et al and Lu et al (\cite{CLWY,LYZ}) have
made some progress in determining which classes of graphs admit
vertex-coloring 2-edge-weighting, notably having shown that
3-connected bipartite graphs are one such class. Khatirineja et al.
\cite{KNN} explored the problem of classifying those graphs which
admit a vertex-coloring 2-edge-weighting. 
Chang et al. \cite{CLWY}  proved that there exists a family of
infinite bipartite graphs (e.g., the generalized $\theta$-graphs)
which are 2-connected and have a vertex-coloring $3$-edge-weighting
but not  vertex-coloring $2$-edge-weightings.

We write  \begin{align*}
 &\mathscr{G}_{12}=\{G\ |\ \mbox{$G$ admits a vertex-coloring
$\{1,2\}$-edge-weighting}\};\\
 &\mathscr{G}_{01}=\{G\ |\ \mbox{$G$ admits a vertex-coloring
$\{0,1\}$-edge-weighting}\};\\
& \mathscr{G}_{12}^{*}=\{G\ |\ \mbox{$G$ is bipartite and admits a
vertex-coloring
 $\{1,2\}$-edge-weighting}\};\\
& \mathscr{G}_{01}^*=\{G\ |\ \mbox{$G$ is bipartite and admits a
vertex-coloring $\{0,1\}$-edge-weighting}\}.
\end{align*}
 Dudek and  Wajc \cite{DW} showed that determining whether a graph belongs to
$\mathscr{G}_{12}$ or $\mathscr{G}_{01}$ is NP-complete. Moreover,
they showed that $\mathscr{G}_{12}\neq \mathscr{G}_{01}$. However,
the counterexamples constructed are non-bipartite. Let $C_6$  be a
cycle of length six and $\Gamma$ be a graph obtained by connecting
an isolated vertex to one of the  vertices of $C_6$. Take two
disjoint copies of $\Gamma$. Connect two vertices of degree one of
two copies and this gives a connected bipartite graph $G$. It is
easy to prove that $G$ admits a vertex-coloring 2-edge-weighting but
not  vertex-coloring $\{0,1\}$-edge-weightings. Hence
$\mathscr{G}_{01}^{*}\neq \mathscr{G}_{12}^{*}$. Next we would like
to propose the following problem.
\begin{problem}
Determining  whether a graph  $G\in \mathscr{G}_{12}^{*}$ or
$G\in\mathscr{G}_{01}^{*}$ is polynomial?
\end{problem}

In this paper,  we characterize bipartite graphs which admit a
vertex-coloring $\mathcal {S}$-edge-weighting  for $\mathcal {S}\in
\{\{0,1\},\{1,2\}\}$, and obtain the following result.
\begin{theorem}\label{02EW}
Let $G$ be a 3-edge-connected bipartite graph $G=(U,W,E)$ with
minimum degree $\delta(G)$. If $G$ contains a vertex $u$ of degree
$\delta(G)$ such that $G-u$ is connected, then  $G$ admits a
vertex-coloring $\mathcal {S}$-edge-weighting for $\mathcal {S}\in
\{\{0,1\},\{1,2\}\}$.
\end{theorem}
By Theorem \ref{02EW}, it is easy to obtain the following result.
\begin{theorem}\label{02EW-conn}
Every 2-connected and 3-edge-connected bipartite graph $G=(U,W,E)$
admits a vertex-coloring $\mathcal {S}$-edge-weighting  for
$\mathcal {S}\in \{\{0,1\},\{1,2\}\}$.
\end{theorem}

So far,  all known counterexamples of bipartite graphs, which do not
have  vertex-coloring $\{0,1\}$-edge-weightings or vertex-coloring
$\{1,2\}$-edge-weightings  are graphs with minimum degree $2$. So we
would like to propose the following problem.
\begin{problem}
Does every bipartite graph with $\delta(G)\geq 3$ admit  a
vertex-coloring $\mathcal {S}$-edge-weighting, where $\mathcal
{S}\in \{\{0,1\},\{1,2\}\}$.
\end{problem}


A \emph{factor} of graph $G$ is a spanning subgraph. For a graph
$G$, there is a close relationship between 2-edge-weighting and
graph factors. Namely, a 2-edge-weighting problem is equivalent to
finding a special factor of graphs (see \cite{ADMRT,ADR}). So to
find factors with pre-specified degree is an important part of
edge-weighting.

Let $g,f:V(G)\rightarrow Z$ be two integer-valued function such that
$g(v)\leq f(v)$ and $g(v)\equiv f(v)\pmod 2$ for all $v\in V(G)$.  A
factor $F$ of $G$ is called \emph{$(g,f)$-parity factor} if
$g(v)\leq d_F(v)\leq f(v)$ and $d_F(v)\equiv f(v)\pmod 2$ for all
$v\in V(G)$. For $X\subseteq V(G)$, we write $g(X)=\sum_{x\in
X}g(x)$ and $f(X)$ is defined similarly. For $(g,f)$-parity factors,
Lov\'asz obtained a sufficient and necessary condition.
\begin{theorem}[Lov\'asz, \cite{Lovasz}]\label{Lovasz}
A graph $G$ contains a $(g,f)$-parity factor if and only if there
exist two disjoint subsets $S$ and $T$ such that
\begin{align*}
\eta(S,T)=f(S)-g(T)+\sum_{x\in T}d_{G-S}(x)-\tau(S,T)<0,
\end{align*}
where $\tau(S,T)$ denotes the number of components $C$, called
$g$-odd components  of $G-S-T$ such that $g(V(C))+e_G(V(C),T)\equiv
1\pmod 2$.
\end{theorem}

In the proof of main theorems, we need the following three Theorems.
\begin{theorem}[Chang et al., \cite{CLWY}]\label{Yu1}
  Every non-trivial connected bipartite graph $G = (A,B, E)$
with $|A|$ even admits a vertex-coloring 2-edge-weighting $w$ such
that $c(u)$ is odd for $u\in A$ and $c(v)$ is even for $v \in B$.
\end{theorem}

\begin{theorem}[Chang et al., \cite{CLWY}]\label{Yu2}
  Let $r\geq 3$ be an integer.  Every $r$-regular bipartite graph
  $G$ admits a vertex-coloring 2-edge-weighting.
\end{theorem}

\begin{theorem}[Khatirinejad et al., \cite{KNN}]\label{KNN}
    Every $r$-regular   graph
  $G$ admits a vertex-coloring 2-edge-weighting
   if and only if it admits a vertex-coloring $\{0,1\}$-edge-weighting.
\end{theorem}

\section{Proof of Theorem \ref{02EW}}

\begin{coro} \label{Yu3}
  Every non-trivial connected bipartite graph $G = (A,B, E)$
with $|A|$ even admits a vertex-coloring $\{0,1\}$-edge-weighting..
\end{coro}

\pf By Theorem \ref{Yu1}, $G$ admits a vertex-coloring
2-edge-weighting $w$ such that $c(u)$ is odd for $u\in A$ and $c(v)$
is even for $v \in B$. Let $w'(e)=0$ if $w(e)=2$ and $w'(e)=1$ if
$w(e)=1$. Then $w'$ is a vertex-coloring $\{0,1\}$-edge-weighting of
graph $G$. \qed

Next, we show that the following two lemmas.

\begin{lemma}\label{lem1}
Let $G$ be a   bipartite graph $G$ with bipartition $(U,W)$, where
$|U|\equiv |W|\equiv 1\pmod 2$. Let $\delta(G)=\delta$ and $u\in U$
such that $d_G(u)=\delta$. If one of two conditions holds, then $G$
contains a factor $F$ such that $d_F(u)=\delta$, $d_F(x)\equiv
\delta+1\pmod 2$ for all $x\in U-u$, $d_F(y)\equiv \delta\pmod 2$
for all $y\in W$ and $d_F(y)\leq \delta-  2$ for all $y\in
N^{\delta}_G(u)$.
\begin{itemize}
\item[$(i)$] $\delta(G)\geq 4$, $G$ is 3-edge-connected and $G-u$ is connected.

\item[$(ii)$] $\delta(G)=3$, $G$ is  3-edge-connected and $|N^{\delta}_G(u)|\leq 2$.

\end{itemize}

%
%
%
%

%
\end{lemma}

\pf Let $M$ be an integer such that $M\geq \triangle(G)$ and
$M\equiv \delta\pmod 2$ and let $m\in \{0,-1\}$ such that $m\equiv
\delta\pmod 2$.   
Let $g, f:V(G)\rightarrow Z$ such that
\[g(x)=\left\{\begin{array}{ll} \delta &\text{if}\
x=u, \\[5pt]
m &\text{if}\ x\in W,  \\[5pt]
m-1 &\text{if} \ x\in U-u,
\end{array}\right.\]
and
\[f(x)=\left\{\begin{array}{ll} M+1 &\text{if}\ x\in U-u,
\\[5pt]
M &\text{if} \ x\in (W\cup \{u\})-N^{\delta}_G(u), \\[5pt]
\delta-2 &\text{if} \ x\in N^{\delta}_G(u).
\end{array}\right.\]

By definition, we have $g(v)\equiv f(v)\pmod 2$ for all $v\in V(G)$.
It is sufficient for us to show that $G$ contains a $(g,f)$-parity
factor. Conversely, suppose that $G$ contains no $(g,f)$-parity
factors. By Theorem \ref{Lovasz}, there exist two disjoint subsets
$S$ and $T$ such that
\begin{align*}
\eta(S,T)=f(S)-g(T)+\sum_{x\in T}d_{G-S}(x)-\tau(S,T)<0,
\end{align*}
where $\tau(S,T)$ denotes the number of   $g$-odd components  of
$G-S-T$. Since $f(V(G))$ is even, by parity, we have
\begin{align}\label{eq:1}
\eta(S,T)=f(S)-g(T)+\sum_{x\in T}d_{G-S}(x)-\tau(S,T)\leq -2.
\end{align}
We choose $S$ and $T$ such that $S\cup T$ is minimal. Let
$A=V(G)-S-T$.

\medskip
\textbf{ Claim 1.~} $S\subseteq N^{\delta}_G(u)$.
\medskip

Otherwise, suppose that $S- N^{\delta}_G(u)\neq \emptyset$ and let
$v\in S-N^{\delta}_G(u)$. Let $S'=S-v$. We have
\begin{align*}
\eta(S',T)&=f(S')-g(T)+\sum_{x\in T}d_{G-S'}(x)-\tau(S',T) \\
&=f(S)-f(v)-g(T)+\sum_{x\in T}d_{G-S}(x)+e_G(v,T)-\tau(S',T)\\
&\leq f(S)-g_1(T)+\sum_{x\in
T}d_{G-S}(x)+e_G(v,T)-(\tau(S,T)-e_G(v,A))-f(v)\\
&=f(S)-g(T)+\sum_{x\in
T}d_{G-S}(x)-\tau(S,T)+(e_G(v,T)+e_G(v,A)-f(v))\\
&\leq f(S)-g(T)+\sum_{x\in
T}d_{G-S}(x)-\tau(S,T)+(d_G(v)-f(v))\\
 &\leq f(S)-g(T)+\sum_{x\in T}d_{G-S}(x)-\tau(S,T)\leq -2,
\end{align*}
contradicting to choice of $S$ and $T$.

\medskip
\textbf{ Claim 2.~} $T\subseteq \{u\}$.
\medskip

Otherwise, let   $v\in T-u$ and $ T'=T-v$. We have
\begin{align*}
\eta(S,T')&=f(S)-g(T')+\sum_{x\in T'}d_{G-S}(x)-\tau(S,T') \\
&=f(S)-g(T)+g(v)+\sum_{x\in T}d_{G-S}(x)-d_{G-S}(v)-\tau(S',T)\\
&\leq f(S)-g(T)+\sum_{x\in
T}d_{G-S}(x)-d_{G-S}(v)-(\tau(S,T)-e_G(v,A))+g(v)\\
&\leq f(S)-g(T)+\sum_{x\in
T}d_{G-S}(x)-\tau(S,T)+(g(v)-d_{G-S}(v)+e_G(v,A))\\
&\leq f(S)-g(T)+\sum_{x\in T}d_{G-S}(x)-\tau(S,T)\leq -2,
\end{align*}
contradicting to choice of $S$ and $T$.

We write $\tau(S,T)=\tau$.  By Claims 1 and 2, (\ref{eq:1}) implies
\begin{align}
\eta(S,T)&=f(S)-g(T)+\sum_{x\in T}d_{G-S}(x)-\tau \\
&=(\delta-2)|S|-\delta |T|+|T|(\delta-|S|)-\tau\\
&=(\delta-2-|T|)|S|-\tau\leq -2\label{eq:2},
\end{align}
which implies $\tau\geq 2$.

 Since $G$ is  3-edge-connected and $G-u$ is connected, then
we have $S\neq \emptyset$ and
\begin{align}\label{eq:3}
3\tau\leq
(\delta-|T|)|S|+|T|(\delta-|S|)=(\delta-2|T|)|S|+|T|\delta.
\end{align}
If $\delta\geq 4$, by (\ref{eq:2}) and (\ref{eq:3}), then  we have
\begin{align}\label{eq:4}
\delta\geq \delta|T|&\geq (2\delta-|T|-6)|S|+6\geq 2\delta-1,
\end{align}
a contradiction if $\delta\geq 4$. So we can assume that $\delta=3$.
Note that  $|S|\leq |N_G^{\delta}(u)|\leq 2$. By (\ref{eq:2}) and
(\ref{eq:3}), we have
\begin{align}\label{eq:41}
3=\delta\geq \delta|T|&\geq -|T||S|+6\geq 4,
\end{align}
a contradiction again.

 This completes the proof. \qed

\begin{lemma}\label{lem2}
Let $G$ be a   bipartite graph $G$ with bipartition $(U,W)$, where
$|U|\equiv |W|\equiv 1\pmod 2$. Let $\delta(G)=\delta$ and $u\in U$
such that $d_G(u)=\delta$. If one of two conditions holds, then $G$
contains a factor $F$ such that $d_F(u)=0$, $d_F(x)\equiv 1\pmod 1$
for all $x\in U-u$, $d_F(y)\equiv 0\pmod 2$   for all $x\in W$ and
$d_F(y)\geq  2$ for all $y\in N_G(u)$.
\begin{itemize}
\item[$(i)$] $\delta(G)\geq 4$, $G$ is 3-edge-connected and $G-u$ is connected.

\item[$(ii)$] $\delta(G)=3$, $G$ is  3-edge-connected and there exists a vertex $v\in N_G(u)$
such that $d_G(v)>3$.

\end{itemize}

\end{lemma}

\pf Let $M$ be an even integer such that $M\geq \triangle(G)$. Let
$g, f:V(G)\rightarrow Z$ such that
\[g(x)=\left\{\begin{array}{ll} 0 &\text{if}\
x\in (\{u\}\cup W)-N_G(u), \\[5pt]
2 &\text{if}\ x\in N_G(u),  \\[5pt]
-1 &\text{if} \ x\in U-u,
\end{array}\right.\]
and
\[f(x)=\left\{\begin{array}{ll} M+1 &\text{if}\ x\in U-u\,
\\[5pt]
M &\text{if} \ x\in W,\\[5pt]
0 &\text{if} \ x=u.
\end{array}\right.\]
Clearly, $g(v)\equiv f(v)\pmod 2$ for all $v\in V(G)$ and $g(V(G))$
is even. It is also sufficient for us to show that $G$ contains a
$(g,f)$-parity factor. Conversely, suppose that $G$ contains no
$(g,f)$-parity factors.   By Theorem \ref{Lovasz}, there exist two
disjoint subsets $S$ and $T$ such that
 \begin{align}\label{eq:5}
\eta(S,T)=f(S)-g(T)+\sum_{x\in T}d_{G-S}(x)-\tau(S,T)\leq -2,
\end{align}
where $\tau(S,T)$ denotes the number of   $g$-odd components  of
$G-S-T$. With similar proof in Lemma \ref{lem1}, we can say
$S\subseteq \{u\}$ and $T\subseteq N_G(u)$. Hence $f(S)=0$ and
$g(T)=2|T|$.  For simplicity, we write $\tau(S,T)=\tau$. By
(\ref{eq:5}), we have
\begin{align}\label{eq:6}
\tau\geq \sum_{x\in T}(d_{G}(x)-|S|)-2|T|+2\geq 2.
\end{align}
Note that $G-u$ is connected, so  we have $|S|\leq 1\leq |T|$. Since
$G$ is 3-edge-connected, then we have
\begin{align}\label{eq:7}
3\tau\leq \sum_{x\in T}(d_{G}(x)-|S|)+(\delta-|T|)|S|.
\end{align}
(\ref{eq:6}) and (\ref{eq:7}) implies
\begin{align}\label{eq:8}
2\sum_{x\in T}d_G(x)+6\leq |S||T|+6|T|+\delta|S|\leq 7|T|+\delta.
\end{align}
If $\delta\geq 4$, then we have $7|T|\geq 6+\delta(2|T|-1)\geq
8|T|+2$, a contradiction. So we can assume that $\delta=3$.
By condition (ii),    $\sum_{x\in T}d_G(x)\geq 3|T|+1$ and $G-u$ is
connected. Combining (\ref{eq:8}),  we have $|T|\geq 5$, a
contradiction since $|T| \leq |N_G(u)|\leq 3$.


 This completes this proof. \qed

\noindent\textbf{Proof of Theorem \ref{02EW}.} By Theorem \ref{Yu1}
and Corollary \ref{Yu3}, we can assume that both $|A|$ and $|B|$ are
odd.

 Firstly, we consider
$\mathcal {S}=\{0,1\}$. If $G$ is 3-regular, by Theorem \ref{Yu2},
then  $G$ admits a vertex-coloring 2-edge-weighting. By Theorem
\ref{KNN}, $G$ also admits a vertex-coloring
$\{0,1\}$-edge-weighting.
%
%
So we can assume that $\delta(G)\geq 3$ and $G$ is not 3-regular. If
$\delta(G)=3$, since $G$ is 3-edge-connected, then $G-x$ is
connected for every vertex $x$ of $G$ with degree no more than four.
Hence   $G$ contains a vertex $v$ with degree three such that $G-v$
is connected and $\sum_{x\in N_G(v)}d_G(v)\geq 3|N_G(v)|+1$. Let
\[u^*=\left\{\begin{array}{ll} u &\text{if}\
\delta \geq 4, \\[5pt]
v &\text{if}\ \delta =3.
\end{array}\right.\]
Without loss generality, we can assume that $u^*\in U$ and so it  is
a vertex satisfying the conditions of Lemma \ref{lem2}. Hence by
Lemma \ref{lem2}, $G$ contains a factor $F$, which satisfies the
following three conditions.
\begin{itemize}
\item[$(i)$]
$d_F(u^*)=0$;

\item[$(ii)$] $d_F(x)\equiv 1\pmod 2$ for all $x\in U-u^*$;

\item[$(iii)$]
$d_F(y)\equiv 0\pmod 2$ for all $x\in W$ and $d_F(y)\geq 2$ for all
$y\in N_G(u^*)$.
\end{itemize}
Clearly, $d_F(x)\neq d_F(y)$ for all $xy \in E(G)$. We assign weight
one for each edge of $E(F)$ and weight zero for each edge of
$E(G)-E(F)$. Then we obtain a vertex-coloring
$\{0,1\}$-edge-weighting of graph $G$.

 Secondly, we show that $G$ admits a vertex-coloring
 2-edge-weighting. By above discussion, we can assume that $G$
 is not 3-regular. If $\delta=3$, since $G$ is 3-edge-connected,
 then $G$ contains a vertex $v'$
 such that  $d_G(v')=3$, $G-v'$ is connected and $|N_G^{\delta}(v')|\leq 2$.

Similarly, we can assume  that $u^*\in U$ is a vertex satisfying the
conditions of Lemma \ref{lem1}. Hence  by Lemma \ref{lem1}, $G$
contains a factor a factor
 $F$ such that
\begin{itemize}
\item[$(i)$] $d_F(u^*)=\delta$;

\item[$(ii)$] $d_F(x)\equiv \delta\pmod 2$ for all $x\in W$
 and $d_F(x)\leq \delta-2$ for  all $x\in N^{\delta}_G(u^*)$;

\item[$(ii)$]  $d_F(y)\equiv \delta+1\pmod 2$ for  all $y\in U-u^*$.
\end{itemize}
Let $w$ be a 2-edge-weighting such that $w(e)=1$ for each $e\in
E(F)$ and $w(e')=2$ for each $e'\in E(G)-E(F)$. Clearly,
$c(u^*)=\delta$. If $y\in N^{\delta}_G(u^*)$, since there exists an
edge $e\sim y$ such that $e\notin E(F)$, then $c(y)=\sum_{e\sim
y}w(e)>\delta$. If $y\in N_G(u^*)- N^{\delta}_G(u^*)$, then
$c(y)\geq d_G(y)>\delta$. Hence $c(y)\neq c(u^*)$ for all $y\in
N_G(u^*)$. For each $xy\in E(G)$, where $x\in U-u^*$ and $y\in W$,
by the choice of $F$, we have $c(x)\equiv \delta+1\pmod 2$ and
$c(y)\equiv \delta\pmod 2$. Hence $w$ is a vertex-coloring
$\{1,2\}$-edge-weighting of graph $G$.

This completes the proof. \qed


\begin{coro}
Let $G$ be a 3-edge-connected bipartite graph. If
$3\leq\delta(G)\leq 5$, then $G$ admits a vertex-coloring $\mathcal
{S}$-edge-weighting for $\mathcal {S}\in \{\{0,1\},\{1,2\}\}$.
\end{coro}

\pf  Since $ 3\leq\delta \leq 5$ and $G$ is 3-edge-connected, then
for every vertex $v$ of degree $\delta$, $G-v$ is connected. By
Lemma \ref{lem1} and Theorem \ref{02EW}, with the same proof, $G$
admits a vertex-coloring $\mathcal {S}$-edge-weighting for $\mathcal
{S}\in \{\{0,1\},\{1,2\}\}$. \qed

\begin{coro}
Let $G$ be a bipartite graph with minimum degree $\delta\geq 4$. If
$G$ is $(\lfloor \frac{\delta}{2}\rfloor+1)$-edge-connected, then
$G$ admits a vertex-coloring $\mathcal {S}$-edge-weighting for
$\mathcal {S}\in \{\{0,1\},\{1,2\}\}$.
\end{coro}




\section{Conclusions}

Let $\mathcal {S}\in \{\{0,1\},\{1,2\}\}$. In this paper, we prove
that every $2$-connected and $3$-edge-connected bipartite graph
admit a vertex-coloring $\mathcal {S}$-edge-weighting. The
generalized $\theta$-graphs  is 2-connected and has a
vertex-coloring $3$-edge-weighting but not a vertex-coloring
$\mathcal {S}$-edge-weightings. So it is an interesting problem to
classify all $2$-connected bipartite graphs admitting a
vertex-coloring $\mathcal {S}$-edge-weighting. Since parity-factor
problem is polynomial, then there exists a polynomial algorithm to
find a vertex-coloring $\mathcal {S}$-edge-weighting of bipartite
graphs satisfying the conditions of  Theorem \ref{02EW}.


\end{document}